\documentclass{amsart}
\usepackage{graphicx}
\usepackage{amssymb}

\newtheorem{conj}{Conjecture}
\newtheorem{cor}{Corollary}
\newtheorem{definition}{Definition}
\newtheorem{identity}{Identity}
\newtheorem{lemma}{Lemma}
\newtheorem{substitution}{Substitution}
\newtheorem{theorem}{Theorem}

\newcommand{\ID}[1] {\vspace{12pt} \begin{identity}#1\end{identity}}
\newcommand{\SUBS}[1] {\vspace{12pt} \begin{substitution}#1\end{substitution}}
\newcommand{\LEMMA}[1] {\vspace{12pt} \begin{lemma}#1\end{lemma}}

\newcommand{\QED}{\hfill $\blacksquare$ \hfill}
\newcommand{\floor}[1] {\left\lfloor #1 \right\rfloor}
\newcommand{\ceil}[1] {\left\lceil #1 \right\rceil}

\title{Boundary slope diameter and crossing number of $2$--bridge knots}
\author{Thomas W.\ Mattman, Gabriel Maybrun \& Kristin Robinson}

\address{Department of Mathematics and Statistics,
         California State University, Chico,
         Chico CA 95929-0525, USA}
\email{TMattman@CSUChico.edu}
\email{GMaybrun@digitalpath.net}
\address{Department of Mathematics,
	 Stanford University,
	  Mathematics, Bldg.\ 380, 
 450 Serra Mall
 Stanford, CA 94305-2125, USA}
\email{krob321@stanford.edu}

\subjclass{Primary 57M25, 57M27}
\keywords{boundary slopes, crossing number, $2$--bridge knots, rational
knots}
\thanks{The second and third authors are undergraduate students who 
worked under the supervision of the first author during an REU at CSU,
Chico in the Summer of 2005. The second author was funded through
NSF REU Award 0354174. The
third author received support from the MAA's program for
Strengthening Underrepresented Minority Mathematics Achievement (SUMMA)
with funding from the NSF, NSA, and Moody's. The first author received
funding through both programs.}

\begin{document}

\begin{abstract}
We prove that for $2$--bridge knots, the diameter, $D$, 
of the set of boundary slopes is twice the crossing number, $c$.  
This constitutes partial verification of a conjecture that, for all knots in $S^3$,
$D \leq 2 c$.
\end{abstract}

\maketitle

\section{Introduction}
Ichihara~\cite{I} told us of a conjecture for knots in $S^3$. Let $D(K)$ denote the 
diameter of the set of boundary slopes of a knot $K$ and $c(K)$ be the crossing number.

\begin{conj} For $K$ a knot in $S^3$, $D(K) \leq 2 c(K)$.
\end{conj}

(To be precise, Ichihara proposed the conjecture only for Montesinos knots and he
and Mizushima~\cite{IM} have recently given a proof of that case.)

Since $0$, being the slope of a Seifert surface, is always included in
the set of boundary slopes, we have, as an immediate consequence, a conjecture due to 
Ishikawa and Shimokawa~\cite{IS}:

\begin{conj} Let $b$ be a finite boundary slope for $K$ a knot in $S^3$. Then $|b| \leq 2 c(K)$.
\end{conj}

For example, it is easy to verify these conjectures for torus knots.
For the unknot, $D(K) = 0 = 2c(K)$. For a non-trivial torus knot $K = (p,q)$ we can assume
$p,q$ relatively prime with $2 \leq q < p$.
The boundary slopes are $0$ and $pq$ \cite{Mo} while the crossing number is
$c(K) = pq-p$~\cite{Mu}.
Thus, $D(K) = pq \leq pq + p(q-2) = 2 c(K)$. Moreover, we have equality for the 
torus $2$--bridge knots which are of the form $(p,2)$ with
$p$ odd. 

We will show that this equality obtains for all $2$--bridge knots:

\begin{theorem} \label{thmain}%
For $K$ a $2$--bridge knot, $D(K) = 2c(K)$.
\end{theorem}

\begin{cor} Let $b$ be a boundary slope for a $2$--bridge knot $K$.
Then $|b| \leq 2 c(K)$.
\end{cor}

This bound is sharp for the $(p,2)$ torus knots and there are many examples
showing that it is also sharp for hyperbolic $2$--bridge knots. 

Using Conway notation, we can associate a rational number $p/q$ to each
$2$--bridge knot $K = K(p/q)$. Hatcher and Thurston~\cite{HT} showed how
to calculate the boundary slopes of $K$ from continued fraction 
representations of $p/q$. On the other hand, the crossing number  
is given by summing the terms in a simple continued fraction for $p/q$
(see~\cite{ES}).
Our technique is, starting with a simple continued fraction for $p/q$, to compute all possible
boundary slope continued fractions and identify those  which yield the maximum and minimum
boundary slopes. We can then verify that the difference
between the maximum and minimum boundary slopes is twice the crossing
number.

In Section~4 we develop four identities for continued fractions
and in Section~5, we use those identities to establish four substitution
rules. These substitution rules will allow us to produce all possible boundary slope
continued fractions for a given rational number:

\begin{theorem} \label{thm1} %
The boundary slope continued fractions of $K( p/q)$ are among the continued fractions obtained
by applying substitutions at non-adjacent positions in the simple continued fraction
of $p/q$. 
\end{theorem}

The proof of Theorem~\ref{thm1} is presented in Section~6 along with the following
corollary.  

\begin{cor} If $\frac p q =[0, a_0, a_1, \ldots, a_n]$ is a simple continued fraction, then $K(
p/q)$ has at most $F_{n+2}$ boundary slopes where $F_n$ is the $n$th Fibonacci number.
\end{cor}

In Section~7 we outline our method for calculating the maximum and
minimum boundary slopes and in Section~8 we prove Theorem~\ref{thmain}.
	
Let us now review the basic ideas of Conway notation and continued
fractions (Section~2) and boundary slopes for $2$--bridge knots (Section~3). 

\section{Conway Notation \& Continued Fractions}

	In this section, we give a brief overview of Conway notation, continued
fractions, and their relationship.
	
 Rational tangles may be constructed by means of tangle algebra (for example, see Adams
\cite{A}).  A rational knot is obtained from numerator closure on a rational tangle.  The
{\em Conway notation} for such a knot is the same as the Conway notation
for the tangle, i.e., a list of integers, $a_0\ a_1\
\dots\ a_n$. Note that the set of rational knots coincides with the set of
$2$--bridge knots and we will use these terms interchangeably.

	A {\em continued fraction} of a rational number $\frac p q$ is a fraction of
the form $$\frac p q = c + \frac{1}{b_0 + \frac{1}{b_1 + \frac{1}{\cdots
+ \frac{1}{b_n}}}} = [c, b_0, b_1, \dots, b_n],$$ where $c \in
\mathbb{Z}$ and each $b_i$, for $0 \leq i \leq n$, is a nonzero integer. 
Note that, since each $b_i$ can be any nonzero integer, the continued
fraction for $\frac p q$ is not unique.  We will call $c$ the {\em integral
component} of the continued fraction, and each $b_i$ will be called a
{\em partial quotient} or {\em term} in the continued fraction.

We will assume that whenever we discuss a continued fraction
$[c, b_0, b_1, \dots, b_n]$, it does, in fact, evaluate to a rational 
number. For example, $[c,2,-1,2]$ is not a valid continued fraction
since
$$c + \frac{1}{2 + \frac{1}{-1+ \frac{1}{2}}} =
c + \frac{1}{2 + \frac{1}{\frac{-1}{2}}} =
c + \frac{1}{2 -2} $$
is not defined as a rational number.

Conway notation and continued fractions are related in that 
we can calculate a rational knot's associated continued fraction by
reversing the order of the Conway notation.  That is, a knot with Conway
notation $a_0\ a_1\ \dots\ a_n$ has the continued fraction $[0, a_n,
a_{n-1}, \dots, a_0] = \frac p q$.  We will denote this knot $K(\frac p
q)$.  In fact (see Cromwell \cite{C} for a proof), all tangles given by
Conway notation corresponding to the same rational number are
equivalent.  Passing to knots introduces additional equivalences:
two rational knots $K(\frac p q)$ and $K(\frac{p'}{q})$
are equivalent if and only if $p' \equiv p^{\pm 1} \pmod q$. Also,
the knots $K( \frac p q)$ and $K( - \frac p q)$ are mirror reflections 
of one another. 
	
Recall that if,  in a
continued fraction $\frac p q = [c, a_0, \dots, a_n]$,
we restrict every $a_i$, $0 \leq i \leq n$, to be a positive
integer, with $a_n > 1$, then this representation of $\frac p q$ is
unique. We will call this the {\em simple continued fraction} of $\frac p q$. 
The corresponding rational knot with Conway notation $a_n\ \dots\ a_0$
then has crossing number $c(K) = \sum_{i=0}^{n}a_i$, as is proven by
Ernst and Sumners \cite{ES}.

	Throughout this paper, we will use the notation $(b_0, \dots, b_m)^c$ 
to mean that 
the pattern ``$b_0, \dots, b_m$'' is repeated $c$ times, with $c$ being
any nonnegative integer, e.g., $[0, (-2, 2)^2] = [0, -2, 2, -2, 2]$ and
$[0, (-2, 2)^0, 2] = [0, 2]$.

		\section{Boundary Slopes}
	
	In this section, we will briefly review how to calculate boundary
slopes for rational knots.

 Let	$B(K)$, or simply $B$, denote
the set of all boundary slopes for a knot $K$.
		For rational knots $K(\frac p q)$, recall \cite{HT} that
$B$ can be calculated from continued
fractions  of $\frac p q$ having every partial quotient at least two in absolute value.
(We will refer to these as {\em boundary slope continued fractions}.)
Specifically, one takes such a continued fraction and pattern-matches
the partial quotients against the pattern
$[+-+-\cdots]$.  The number of terms matching this pattern we call $b^+$,
and the number of terms not matching this pattern (e.g., the total number
of terms minus $b^+$) we call $b^-$ (since these terms match the pattern
$[-+-+\cdots])$. In this way, we associate to each boundary slope continued fraction
two non-negative integers $b^+$ and $b^-$. 

Among the boundary slope continued fractions, there is a unique one consisting
only of even terms (indeed, this is easy to see using the four
substitutions we will derive in Section 5). This corresponds
to a Seifert surface of boundary slope $0$, so we will denote its $b^+$ and $b^-$ by $b_0^+$ and
$b_0^-$. Then, the boundary slope associated to any other continued fraction is
given by comparing its $b^+$ and $b^-$ with those of the Seifert slope;
the boundary slope is $  2\big(
(b^+-b^-) - (b_0^+ - b_0^-)\big)$. Applying this calculation to
every continued fraction with terms at least two in absolute value gives
the set of boundary slopes $B$. $B$ is a finite set of even integers.
The {\em diameter} $D(K)$ is the difference between the maximum and minimum elements
of $B$.

\section{Continued Fraction Identities}

	In this section, we will prove four identities related to continued
fractions.  
	For identities 2 and 4, we will allow the last entry in a continued
fraction to be any nonzero rational number, provided the resulting
continued fraction represents a rational number.  Note
that $[b_0, \dots, b_m, k] = [b_0, \dots, b_m, a_0, \dots,
a_n]$ when $k = [a_0, \dots, a_n]$.
	
	Throughout this section, let $\mathbb{N}_0 = \mathbb{N} \cup \{0\}$ and
$\mathbb{Q}^* = \mathbb{Q} \setminus \{0\}$.
	
\begin{identity}
 Let $c \in \mathbb{N}$.  Then $$[(-2, 2)^c] = -\frac{2c+1}{2c}$$ 
\end{identity}

	\proof  By induction.
	
	\textbf{Base Case ($c=1$):}  $[-2, 2] = -2 + \frac 1 2 = -\frac 3 2 = -\frac{2 \cdot 1 + 1}{2 \cdot 1}$.
	
	\textbf{Induction Step:}  Assume that $[(-2, 2)^c] = \frac{2c+1}{2c}$.  Then
	\begin{eqnarray*}
	 [(-2, 2)^{c+1}] 
		&=& -2 + \frac{1}{2 + \frac{1}{[(-2, 2)^c]}} \\
		&=& -2 + \frac{1}{2 + \frac{1}{-\frac{2c+1}{2c}}} = -2 +
\frac{2c+1}{2c+2} \\
		&=& -\frac{2c+3}{2c+2}		= -\frac{2(c+1)+1}{2(c+1)}
	\end{eqnarray*} \QED
	
	\ID{ Let $c \in \mathbb{N}_0$ and $k \in \mathbb{Q}^*$.  Then $$[(-2, 2)^c, k] = \frac{2ck +2c
+k}{1 - 2ck - 2c}$$ }
	
Note that the denominator becomes zero only in the case where the
continued fraction does not represent a rational number.

	\proof By induction.
	
	\textbf{Base Case ($c=0$):}  $[k] = k = \frac{2 \cdot 0 \cdot k + 2 \cdot 0 + k}{1 - 2 \cdot 0 \cdot k - 2 \cdot 0}$.
	
	\textbf{Induction Step:}  Assume that $[(-2, 2)^c, k] = \frac{2ck +2c +k}{1 - 2ck - 2c}$.  Then
	\begin{eqnarray*}
		[(-2, 2)^{c+1}, k] 
		&=& -2 + \frac{1}{2 + \frac{1}{[(-2, 2)^c, k]}} \\
		&=& -2 + \frac{1}{2 + \frac{1}{\frac{2ck +2c +k}{1 - 2ck - 2c}}} \\
		&=& -2 + \frac{2ck + 2c + k}{2ck + 2c + 2k + 1} \\
		&=& \frac{-2ck - 2c - 3k - 2}{2ck + 2c + 2k + 1} \\
		&=& \frac{2(c+1)k + 2(c+1) + k}{1 - 2(c+1)k - 2(c+1)}
	\end{eqnarray*} \QED
	
	\ID{ Let $c \in \mathbb{N}$.  Then $$[(2, -2)^c] = \frac{2c+1}{2c}$$ }

	\proof  Note that $[(2, -2)^c] = [2, (-2, 2)^{c-1}, -2]$.  Apply Identity 2.  \QED
	
	\ID{ Let $c \in \mathbb{N}_0$ and $k \in \mathbb{Q}^*$.  Then $$[(2, -2)^c, k] = \frac{2ck - 2c
+ k}{2ck - 2c + 1}$$ }

Again, the denominator becomes zero only in the case where the
continued fraction does not represent a rational number.

	\proof  This proof will be done in two parts.

	\textbf{Case 1 ($c=0$):}  $[k] = k = \frac{2 \cdot 0 \cdot k - 2 \cdot 0
+ k}{2 \cdot 0 \cdot k - 2 \cdot 0 + 1}$.
	
	\textbf{Case 2 ($c>0$):}  Note that $[(2, -2)^c, k] = [2, (-2,
2)^{c-1}, -2, k]$.  Apply Identity 2.  \QED

\section{Substitution Rules}

	In this section, we will prove four identities, or substitutions,
which will be used to derive equal continued fractions.  In particular,
given $\frac p q = [c, a_0, \dots, a_n]$, where every $a_i$ is
positive, these substitutions can be used to
calculate all the boundary slope continued fractions of $\frac p q$, i.e., continued fractions
$[c', b_0, \dots, b_m]$ where each $|b_i| \geq 2$.
We conclude the section with an example to illustrate how these rules can
be applied to a specific continued fraction.

	Throughout this section, let $\mathbb{N}_0 = \mathbb{N} \cup \{0\}$ and
$\mathbb{Z}^* = \mathbb{Z} \setminus \{0\}$.
	
	\SUBS { Let $n \in \mathbb{N}$.  Let $a_0 \in \mathbb{Z}$ and $a_1 \in \mathbb{N}$.  If $n = 2$ then let $a_2 \in \mathbb{Z} \setminus \{0, -1\}$.  If $n \geq 3$ then let $a_i \in \mathbb{Z}^*$ for all $2 \leq i \leq n$.  If we have $r = [a_0, 2a_1, a_2, a_3, \dots, a_n]$, then $r = [a_0+1, (-2, 2)^{a_1-1}, -2, a_2+1, a_3, a_4, \dots, a_n]$. }
	\proof  This proof will be done in three parts.
	
	\textbf{Case 1 ($n = 1$):}  We want to show that $[a_0, 2a_1] = [a_0+1, (-2, 2)^{a_1-1}, -2]$.
	\begin{eqnarray*}
	 && [a_0+1, (-2, 2)^{a_1-1}, -2] \\
		&=& a_0+1+\frac{1}{[(-2, 2)^{a_1-1}, -2]} \\
		&=& a_0+1+\frac{1}{\frac{2(a_1-1)(-2) + 2(a_1-1) +
(-2)}{1-2(a_1-1)(-2)-2(a_1-1)}} \mbox{\qquad (Apply Identity 2)} \\
		&=& a_0+1+\frac{-2a_1+1}{2a_1} \\
		&=& a_0 + \frac{1}{2a_1} \\
		&=& [a_0, 2a_1]
	\end{eqnarray*}


	\textbf{Case 2 ($n = 2$):}  We want to show $[a_0, 2a_1, a_2] =
[a_0+1, (-2, 2)^{a_1-1}, -2, a_2+1]$.
	\begin{eqnarray*}
		&& [a_0+1, (-2, 2)^{a_1-1}, -2, a_2+1] \\
		&=& a_0+1+\frac{1}{[(-2, 2)^{a_1-1}, -2, a_2+1]} \\
		&=& a_0+1+\frac{1}{\frac{2(a_1-1)\left(\frac{-2a_2-1}{a_2+1}\right) + 2(a_1-1) + \frac{-2a_2-1}{a_2+1}}{1-2(a_1-1)\left(\frac{-2a_2-1}{a_2+1}\right)-2(a_1-1)}} \mbox{\qquad (Apply Identity 2)} \\
		&=& a_0+1+\frac{a_2-2a_1a_2-1}{2a_1a_2+1} \\
		&=& a_0 + \frac{1}{\frac{2a_1a_2+1}{a_2}} \\
		&=& a_0 + \frac{1}{2a_1 + \frac{1}{a_2}} \\
		&=& [a_0, 2a_1, a_2]
	\end{eqnarray*}

	\textbf{Case 3 ($n \geq 3$):}  We want to show that $[a_0, 2a_1, a_2, a_3, a_4, \dots, a_n] = [a_0+1, (-2, 2)^{a_1-1}, -2, a_2+1, a_3, a_4, \dots, a_n]$.  Let $R = [a_3, a_4, \dots, a_n]$.
	\begin{eqnarray*}
		&& [a_0+1, (-2, 2)^{a_1-1}, -2, a_2+1, a_3, a_4, \dots, a_n] \\
		&=& a_0+1 + \frac{1}{[(-2, 2)^{a_1-1}, -2, a_2+1, a_3, a_4, \dots, a_n]} \\
		&=& a_0+1 + \frac{1}{\frac{2(a_1-1)\left(\frac{-2Ra_2-R-2}{Ra_2+R+1}\right) + 2(a_1-1) + \frac{-2Ra_2-R-2}{Ra_2+R+1}}{1 - 2(a_1-1)\left(\frac{-2Ra_2-R-2}{Ra_2+R+1}\right) - 2(a_1-1)}}  \mbox{\qquad (Apply Identity 2)} \\
		&=& a_0+1 + \frac{Ra_2 - 2a_1 - R - 2Ra_1a_2 + 1}{R + 2a_1 + 2Ra_1a_2} \\
		&=& a_0 + \frac{1}{\frac{R+2a_1+2Ra_1a_2}{Ra_2+1}} \\
		&=& a_0 + \frac{1}{2a_1 + \frac{R}{Ra_2+1}} \\
		&=& a_0 + \frac{1}{2a_1 + \frac{1}{a_2 + \frac{1}{R}}} \\
		&=& [a_0, 2a_1, a_2, a_3, a_4, \dots, a_n]
\end{eqnarray*}

\QED
	
	\SUBS { Let $n \in \mathbb{N}$.  Let $a_0 \in \mathbb{Z}$ and $a_1 \in \mathbb{N}$.  If $n = 2$ then let $a_2 \in \mathbb{Z} \setminus \{0, 1\}$.  If $n \geq 3$ then let $a_i \in \mathbb{Z}^*$ for all $2 \leq i \leq n$.  If we have $r = [a_0, -2a_1, a_2, a_3, \dots, a_n]$, then $r = [a_0-1, (2, -2)^{a_1-1}, 2, a_2-1, a_3, a_4, \dots, a_n]$. }
	\proof  This proof will be done in three parts.

	\textbf{Case 1 ($n = 1$):}  We want to show that $[a_0, -2a_1] = [a_0-1, (2, -2)^{a_1-1}, 2]$.
	\begin{eqnarray*}
		 [a_0-1, (2, -2)^{a_1-1}, 2] 
		&=& a_0-1+\frac{1}{[(2, -2)^{a_1-1}, 2]} \\
		&=& a_0-1+\frac{1}{\frac{2(a_1-1)(2) - 2(a_1-1) + (2)}{2(a_1-1)(2)-2(a_1-1) + 1}} \mbox{\quad
(Apply Identity 4)} \\
		&=& a_0-1+\frac{2a_1-1}{2a_1} \\
		&=& a_0 + \frac{1}{-2a_1} \\
		&=& [a_0, -2a_1]
	\end{eqnarray*}


	\textbf{Case 2 ($n = 2$):}  We want to show $[a_0, -2a_1, a_2] = [a_0-1, (2, -2)^{a_1-1},
2, a_2-1]$.
	\begin{eqnarray*}
		&& [a_0-1, (2, -2)^{a_1-1}, 2, a_2-1] \\
		&=& a_0-1+\frac{1}{[(2, -2)^{a_1-1}, 2, a_2-1]} \\
		&=& a_0-1+\frac{1}{\frac{2(a_1-1)\left(\frac{2a_2-1}{a_2-1}\right) - 2(a_1-1) + \frac{2a_2-1}{a_2-1}}{2(a_1-1)\left(\frac{2a_2-1}{a_2-1}\right)-2(a_1-1) + 1}} \mbox{\qquad (Apply Identity 4)} \\
		&=& a_0-1+\frac{2a_1a_2-a_2-1}{2a_1a_2-1} \\
		&=& a_0 + \frac{1}{\frac{2a_1a_2-1}{-a_2}} \\
		&=& a_0 + \frac{1}{{-2a_1}+\frac{1}{a_2}} \\
		&=& [a_0, -2a_1, a_2]
	\end{eqnarray*}

	\textbf{Case 3 ($n \geq 3$):}  We want to show that $[a_0, -2a_1, a_2, a_3, a_4, \dots, a_n] = [a_0-1, (2, -2)^{a_1-1}, 2, a_2-1, a_3, a_4, \dots, a_n]$.  Let $R = [a_3, a_4, \dots, a_n]$.
	\begin{eqnarray*}
		&& [a_0-1, (2, -2)^{a_1-1}, 2, a_2-1, a_3, a_4, \dots, a_n] \\
		&=& a_0-1
+ \frac{1}{[(2, -2)^{a_1-1}, 2, a_2-1, a_3, a_4, \dots, a_n]} \\
		&=& a_0-1 +
\frac{1}{\frac{2(a_1-1)\left(\frac{2Ra_2+2-R}{Ra_2+1-R}\right) - 2(a_1-1)
+
\frac{2Ra_2+2-R}{Ra_2+1-R}}{2(a_1-1)\left(\frac{2Ra_2+2-R}{Ra_2+1-R}\right)
- 2(a_1-1)+1}}  \mbox{\qquad (Apply Identity 4)} \\
		&=& a_0-1 + \frac{2a_1-R-Ra_2+2Ra_1a_2-1}{2a_1-R+2Ra_1a_2} \\
		&=& a_0 + \frac{1}{\frac{R-2a_1-2Ra_1a_2}{Ra_2+1}} \\
		&=& a_0 + \frac{1}{-2a_1 + \frac{R}{Ra_2+1}} \\
		&=& a_0 + \frac{1}{-2a_1 + \frac{1}{a_2 + \frac{1}{R}}} \\
		&=& [a_0, -2a_1, a_2, a_3, a_4, \dots, a_n]
	\end{eqnarray*}

\vspace{-12pt}	

\QED

	\SUBS { Let $n \in \mathbb{N}$.  Let $a_0 \in \mathbb{Z}$ and $a_1 \in \mathbb{N}_0$.  If $n =
2$ then let $a_2 \in \mathbb{Z} \setminus \{0, -1\}$.  If $n \geq 3$ then let $a_i \in
\mathbb{Z}^*$ for all $2 \leq i \leq n$.  If we have $r = [a_0, 2a_1+1, a_2, a_3, \dots, a_n]$,
then $r = [a_0+1, (-2, 2)^{a_1}, -a_2-1, -a_3, -a_4, \dots, -a_n]$. }
	\proof  This proof will be done in three parts.
	
	\textbf{Case 1 ($n = 1$):}  We want to show that $[a_0, 2a_1+1] = [a_0+1, (-2, 2)^{a_1}]$.
	Note: When $a_1=0$ this is trivially true.  So, we can assume $a_1>0$.
	\begin{eqnarray*}
	 [a_0+1, (-2, 2)^{a_1}] 
		&=& a_0+1+\frac{1}{[(-2, 2)^{a_1}]} \\
		&=& a_0+1+\frac{1}{-\frac{2a_1+1}{2a_1}} \mbox{\qquad (Apply Identity 1)} \\
		&=& a_0+1+\frac{2a_1}{-2a_1-1} \\
		&=& a_0 + \frac{1}{2a_1+1} \\
		&=& [a_0, 2a_1+1]
	\end{eqnarray*}
	

	\textbf{Case 2 ($n = 2$):}  We want to show $[a_0, 2a_1+1, a_2] = [a_0+1, (-2, 2)^{a_1},-
a_2-1]$.
	\begin{eqnarray*}	
	&& [a_0+1, (-2, 2)^{a_1}, -a_2-1]  \\
		&=& a_0+1+\frac{1}{[(-2, 2)^{a_1}, -a_2-1]} \\
		&=& a_0+1+\frac{1}{\frac{2a_1(-a_2-1)+2a_1+(-a_2-1)}{1-2a_1(-a_2-1)-2a_1}} \mbox{\qquad (Apply
Identity 2)} \\
		&=& a_0+1+\frac{-2a_1a_2-1}{2a_1a_2+a_2+1} \\
		&=& a_0 + \frac{1}{\frac{2a_1a_2+a_2+1}{a_2}} \\
		&=& a_0 + \frac{1}{2a_1+1+\frac{1}{a_2}}\\
		&=& [a_0, 2a_1+1, a_2]
	\end{eqnarray*}
	
	\textbf{Case 3 ($n \geq 3$):}  We want to show that $[a_0, 2a_1+1, a_2, a_3, a_4, \dots, a_n] = [a_0+1, (-2, 2)^{a_1}, -a_2-1, -a_3, -a_4, \dots, -a_n]$.  Let $R = [a_3, a_4, \dots, a_n]$.  Note: $-R = [-a_3, -a_4, \dots, -a_n]$.
	\begin{eqnarray*}
		&& [a_0+1, (-2, 2)^{a_1}, -a_2-1, -a_3, -a_4, \dots, -a_n] \\
		&=& a_0+1 + \frac{1}{[(-2, 2)^{a_1}, -a_2-1, -a_3, -a_4, \dots, -a_n]} \\
		&=& a_0+1 + \frac{1}{\frac{2a_1\left(\frac{Ra_2+R+1}{-R}\right) + 2a_1 + \frac{Ra_2+R+1}{-R}}{1-2a_1\left(\frac{Ra_2-R+1}{-R}\right)-2a_1}}  \mbox{\qquad (Apply Identity 2)} \\
		&=& a_0+1 + \frac{-R-2a_1-2Ra_1a_2}{2a_1+R+Ra_2+2Ra_1a_2+1} \\
		&=& a_0 + \frac{1}{\frac{2a_1+R+Ra_2+2Ra_1a_2+1}{1+Ra_2}} \\
		&=& a_0 + \frac{1}{2a_1+1 + \frac{1}{\frac{Ra_2+1}{R}}} \\
		&=& a_0 + \frac{1}{2a_1+1+\frac{1}{a_2+\frac{1}{R}}} \\
		&=& [a_0, 2a_1+1, a_2, a_3, a_4, \dots, a_n] 
	\end{eqnarray*} 
\QED
	
	\SUBS { Let $n \in \mathbb{N}$.  Let $a_0 \in \mathbb{Z}$ and $a_1 \in \mathbb{N}_0$.  If $n =
2$ then let $a_2 \in \mathbb{Z}\setminus\{0,1\}$.  If $n \geq 3$ then let $a_i \in \mathbb{Z}^*$
for all $2 \leq i \leq n$.  If we have $r = [a_0, -2a_1-1, a_2, a_3, \dots, a_n]$, then $r =
[a_0-1, (2, -2)^{a_1}, -a_2+1, -a_3, -a_4, \dots, -a_n]$. }
	\proof  This proof will be done in three parts.
	
	\textbf{Case 1 ($n = 1$):}  We want to show that $[a_0, -2a_1-1] = [a_0-1, (2, -2)^{a_1}]$.
	Note: When $a_1=0$ this is trivially true.  So, we can assume $a_1>0$.
	\begin{eqnarray*}
		 [a_0-1, (2, -2)^{a_1}] 
		&=& a_0-1+\frac{1}{[(2, -2)^{a_1}]} \\
		&=& a_0-1+\frac{1}{\frac{2a_1+1}{2a_1}} \mbox{\qquad (Apply Identity 3)} \\
		&=& a_0-1+\frac{2a_1}{2a_1+1} \\
		&=& a_0 + \frac{1}{-2a_1-1} \\
		&=& [a_0, -2a_1-1]
	\end{eqnarray*}
	

	\textbf{Case 2 ($n = 2$):}  We want to show $[a_0, -2a_1-1, a_2] = [a_0-1, (2, -2)^{a_1},-
a_2+1]$.
	\begin{eqnarray*}
		&& [a_0-1, (2, -2)^{a_1}, -a_2+1] \\
		&=& a_0-1+\frac{1}{[(-2, 2)^{a_1}, -a_2+1]} \\
		&=& a_0-1+\frac{1}{\frac{2a_1(-a_2+1)-2a_1+(-a_2+1)}{2a_1(-a_2+1)-2a_1+1}} \mbox{\qquad (Apply Identity 4)} \\
		&=& a_0-1+\frac{2a_1a_2-1}{2a_1a_2+a_2-1} \\
		&=& a_0 + \frac{1}{\frac{1-a_2-2a_1a_2}{a_2}} \\
		&=& a_0 + \frac{1}{-2a_1-1+\frac{1}{a_2}}\\
		&=& [a_0, -2a_1-1, a_2]
	\end{eqnarray*}
	
	\textbf{Case 3 ($n \geq 3$):}  We want to show that $[a_0, -2a_1-1, a_2, a_3, a_4, \dots, a_n] = [a_0-1, (2, -2)^{a_1}, -a_2+1, -a_3, -a_4, \dots, -a_n]$.  Let $R = [a_3, a_4, \dots, a_n]$.  Note: $-R = [-a_3, -a_4, \dots, -a_n]$.
	\begin{eqnarray*}
		&& [a_0-1, (2, -2)^{a_1}, -a_2+1, -a_3, -a_4, \dots, -a_n] \\
		&=& a_0-1 + \frac{1}{[(2, -2)^{a_1}, -a_2+1, -a_3, -a_4, \dots, -a_n]} \\
		&=& a_0-1 + \frac{1}{\frac{2a_1\left(\frac{R-Ra_2-1}{R}\right) - 2a_1 + \frac{R-Ra_2-1}{R}}{2a_1\left(\frac{R-Ra_2-1}{R}\right) - 2a_1 + 1}}  \mbox{\qquad (Apply Identity 4)} \\
		&=& a_0-1 + \frac{2a_1-R+2Ra_1a_2}{2a_1-R+Ra_2+2Ra_1a_2+1} \\
		&=& a_0 + \frac{1}{\frac{R-2a_1-Ra_2-2Ra_1a_2-1}{Ra_2+1}} \\
		&=& a_0 + \frac{1}{-2a_1-1 + \frac{1}{\frac{Ra_2+1}{R}}} \\
		&=& a_0 + \frac{1}{-2a_1-1+\frac{1}{a_2+\frac{1}{R}}} \\
		&=& [a_0, -2a_1-1, a_2, a_3, a_4, \dots, a_n]
	\end{eqnarray*} \QED

%
%
%
		
	\subsection{An example of the application of Substitutions 1--4}
	
Let us illustrate how the above results can be used to generate a list
of all boundary slope continued fractions starting from
the simple continued fraction. As an example, suppose
we start with
$[0, 2a, 2b+1, 2c]$, where $a,c \in \mathbb{N}$ and $b \in \mathbb{N} \cup \{0\}$.  By
applying Substitution 1, we can immediately derive another continued fraction: $[1, (-2,
2)^{a-1}, -2, 2b+2, 2c]$.  We will refer to this as {\em applying Substitution 1 at position 0}
as it is the $a_0$ term, $2a$, that has been replaced by the sequence $-2, 2, \dots, -2$. 

Applying the same substitution at position 2, we get $[1, (-2,
2)^{a-1}, -2, 2b+3, (-2, 2)^{c-1}, -2]$. We could continue on this path, but it is easy to see
that any further substitutions will result in a $\pm 1$ term.  Therefore, we return to the
original sequence and use Substitution 3 (at position 1) to obtain $[0, 2a+1, (-2, 2)^b,
-2c-1]$.  Finally, applying Substitution 1 at position 2, we have 
$[0, 2a, 2b+2, (-2, 2)^{c-1}, -2]$.

Thus, there are five boundary slope continued fractions that can be derived from
the simple continued fraction $[0, 2a, 2b+1, 2c]$: three obtained by substitutions at
positions 0, 1, and  2; one by substitutions at 0 and 2; and the original continued fraction 
itself (with no substitutions). Note, that these are precisely the fractions obtained by applying
substitutions at non-adjacent positions.
	
\section{Proof of Theorem~\ref{thm1}}

In this section we will prove Theorem~\ref{thm1}, that 
the boundary slope continued fractions are among the fractions
obtained by applying substitutions at non-adjacent positions in the original simple continued
fraction. Our strategy is to first review Langford's argument~\cite{L}
that the boundary slopes are determined by the leaves of a binary tree. 
We then show, by induction, that applying substitutions at non-adjacent positions 
accounts for all the leaves of the tree.

\subsection{The boundary slope binary tree}

Recall that for any rational $\frac
p q$, we can find another rational $\frac{p'}{q}$ such that $0 \leq \frac{p'}{q} < 1$ and
$K\left(\frac p q\right) = K\left(\frac{p'}{q}\right)$.  Also, recall that the rational $0$
corresponds to the unknot, which has a rather boring set of continued fractions (namely, $0$ is
the only one).  Therefore, without loss of generality, we will assume henceforth that $0 < \frac
p q < 1$.  There is a unique simple continued fraction, $[0, a_0, \dots, a_n]$, for such
$\frac p q$, such that $a_n \geq 2$ and, for all $i \in \{0, 1, \dots, n\}$, $a_i > 0$.
	
Before we can prove Theorem~\ref{thm1}, we must first state a lemma.
The straightforward proof by induction may be found in Langford~\cite{L} which is
also the source for the following definition.

\begin{definition}	
The {\em $k$th subexpansion} of $[c, a_0, \dots,
a_n]$ is the continued fraction $[0, a_k, \dots, a_n]$ where $0 \leq k \leq n$.  
\end{definition}

	\LEMMA{Let $[c, a_0, \dots, a_n]$ be a boundary slope continued fraction, that is, for
each $i \in \{0, 1,
\dots, n\}$, $\left|a_i\right| \geq 2$.  Then every subexpansion $r$ of $[c, a_0, \dots, a_n]$
satisfies $\left|r\right| < 1$.}

As Langford \cite{L} has shown, a complete list of boundary slope continued fractions for a
rational $\frac p q$, where each partial quotient is at least two in absolute value, can be
calculated by means of a binary tree.  We will now outline the creation of this binary tree
which follows from Lemma~1. 

The root vertex is labelled with the fraction $\frac p q$ and the two edges coming from
the root are labelled $0 = \lfloor \frac p q \rfloor$ and $1 = \lceil \frac p q \rceil$.
At every other vertex in the tree, we arrive with the first $k$ terms in a continued fraction for
$\frac p q$ and a rational number $r$ representing the $(k-1)$st subexpansion. The $k$ terms are
found as labels of the edges of the tree starting from the root and continuing to the vertex in
question. We label the vertex with $r$. Since, by Lemma~1, any $k$th subexpansion is less than
one in absolute value, we know that the next term in the continued fraction, $a_i$, is within
$1$ of
$1/r$: $ |a_i-1/r| < 1$. However, $a_i$ is an integer. Therefore, $a_i$ is either the
floor 
$\lfloor 1/r \rfloor$ or the ceiling $\lceil 1/r \rceil$ of $1/r$. If $1/r$ is not an
integer, there will be two edges coming out of the vertex, one labelled with $\lfloor 1/r
\rfloor$, and the other labelled with $\lceil 1/r \rceil$. Since $|r| < 1$, neither of these
arrows is $0$.  If either is $\pm 1$,
we terminate that edge with a leaf labelled ``$\nexists$" to
indicate that this path does not lead to a  boundary slope continued fraction.
(When we refer to the leaves of the binary tree below, we will be excluding
these ``dead" leaves.) If $1/r$ is an
integer, then, there is only one edge coming out of the vertex. Label the edge with $1/r$
and label the leaf vertex at the end of this edge with the continued fraction expansion for
$\frac p q$ given by the labels of the edges from the root to the leaf. 

\begin{figure}[h]
\begin{center}
\includegraphics[scale = 0.28]{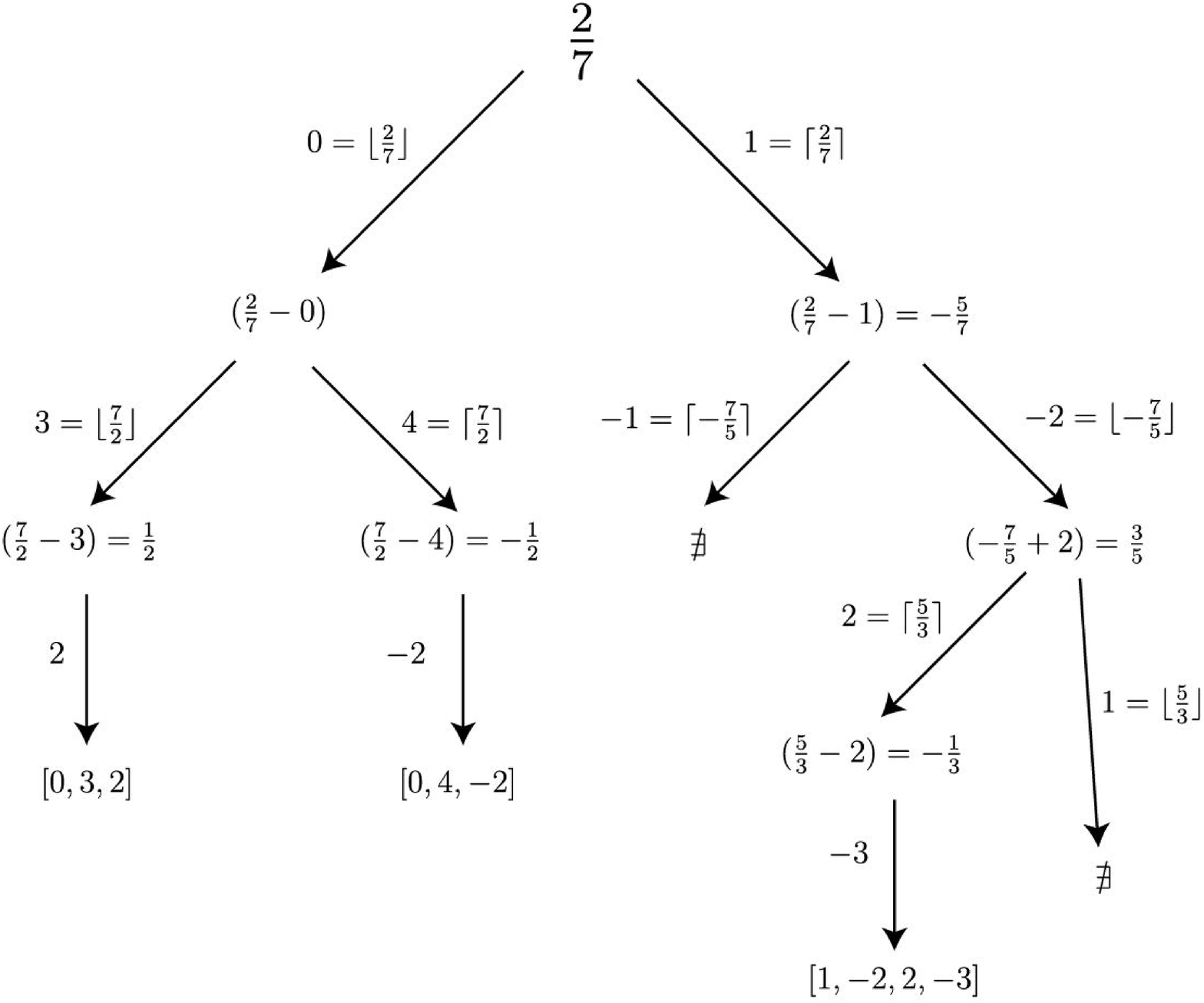}
\caption{\label{figbtree}%
The boundary slope binary tree for $\frac p q = \frac 2 7$ (the $5_2$ knot).}
\end{center}
\end{figure}

For example, Figure~\ref{figbtree}
shows the binary tree for the fraction $2/7$ (which corresponds to the $5_2$ knot).

Thus, by Lemma~1, the algorithm
used to construct the tree will provide all the boundary slope continued fractions of
$\frac p q $ as leaf vertices.

\subsection{Binary tree from substitutions}
Now, let's prove the theorem by 
showing that the leaves of  Langford's binary tree (and therefore the set of boundary
slopes) correspond to applying substitutions  at non-adjacent positions in the simple
continued fraction.

\setcounter{theorem}{1}

\begin{theorem}
The boundary slope continued fractions of $K( p / q )$ are among the continued fractions
obtained by applying substitutions at non-adjacent positions in the simple continued fraction
of $p/q $. 
\end{theorem}

\proof
We proceed by induction on the length $n$ of the simple continued fraction
$[0, a_0, a_1, \ldots, a_n]$.

\textbf{Case 1} ($n=0$):
Here, $p/q = 1/a_0$. We wish to show that the boundary slope continued fractions are
among the two continued fractions given by substituting or not at position 0. There are three
subcases. (To simplify the exposition, we will not be considering the, very similar, trees 
that arise when the terms $a_i$ are negative
although they may be required as part of our induction.) 

\begin{figure}[h]
\begin{center}
\includegraphics[scale=0.09]{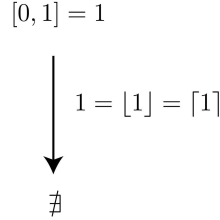}
\caption{\label{fig1}%
The binary tree for $[0,1]$.}
\end{center}
\end{figure}

\underline{Subcase 1} ($a_0 = 1$): In this case, the tree is shown in Figure~\ref{fig1}. 
There are no boundary slope continued fractions in this case.
(Actually, here $\frac p q = 1$, so we've violated our assumption that 
$\frac p q < 1$. Ordinarily, we would represent this knot, the unknot, by
$[0]$ and that would also be the only boundary slope.
We include this case as it may arise as part of our induction.)
Thus, it is true that all boundary slope
continued fractions are among the two continued fractions $[0,1]$ and $[1]$ given by
substituting or not at position 0.

\begin{figure}[h]
\begin{center}
\includegraphics[scale = 0.43]{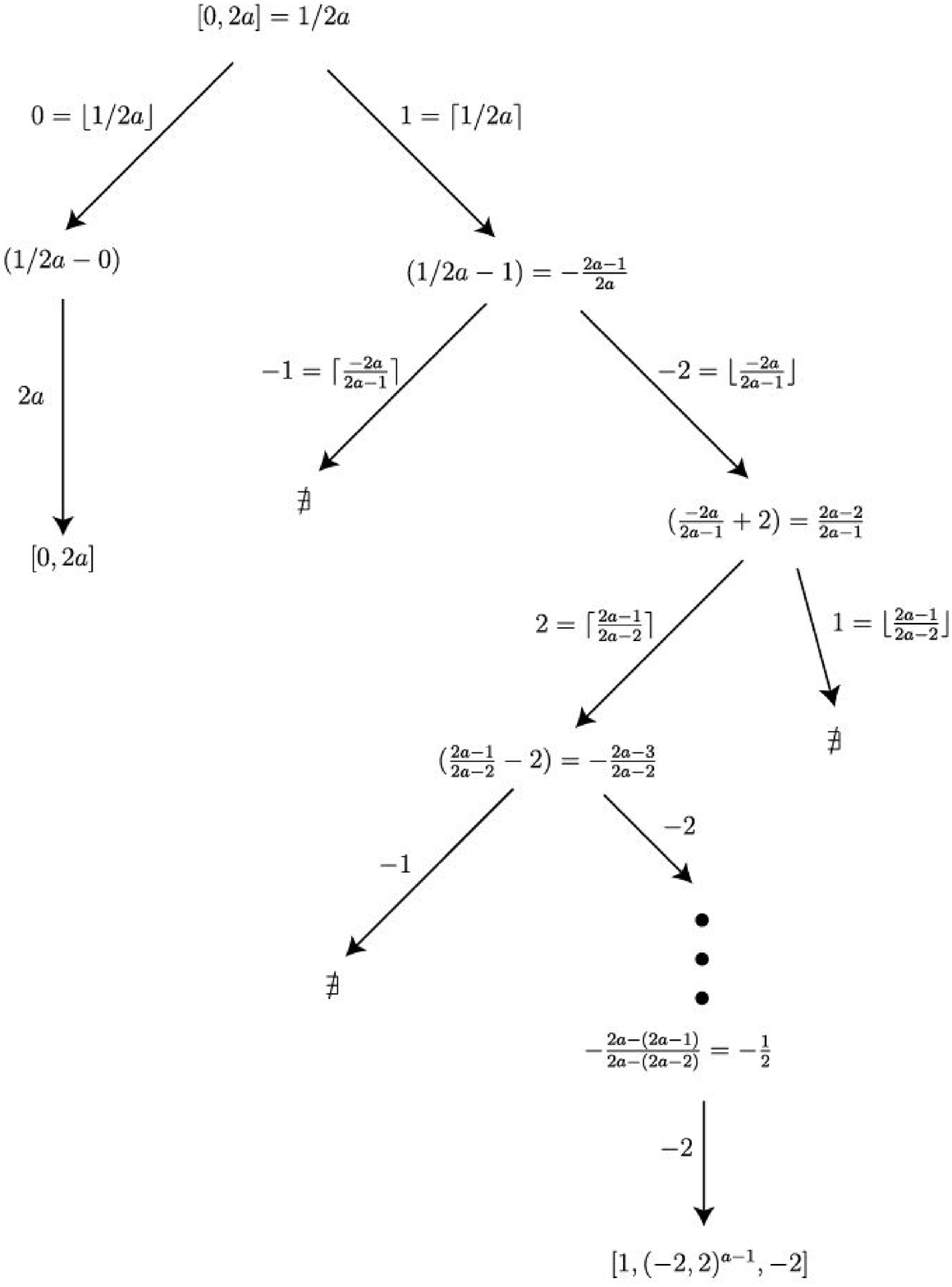}
\caption{\label{fig2}%
The binary tree for $[0,2a]$.}
\end{center}
\end{figure}

\underline{Subcase 2} ($a_0 = 2a$, $a \geq 1$): 
The binary tree is shown in Figure~\ref{fig2}. 
There are two
boundary slope continued fractions, and they are the fractions $[0, a_0]$ and
$[1, (-2,2)^a, -2]$ given by substituting or not at position 0.

\begin{figure}[h]
\begin{center}
\includegraphics[scale=0.39]{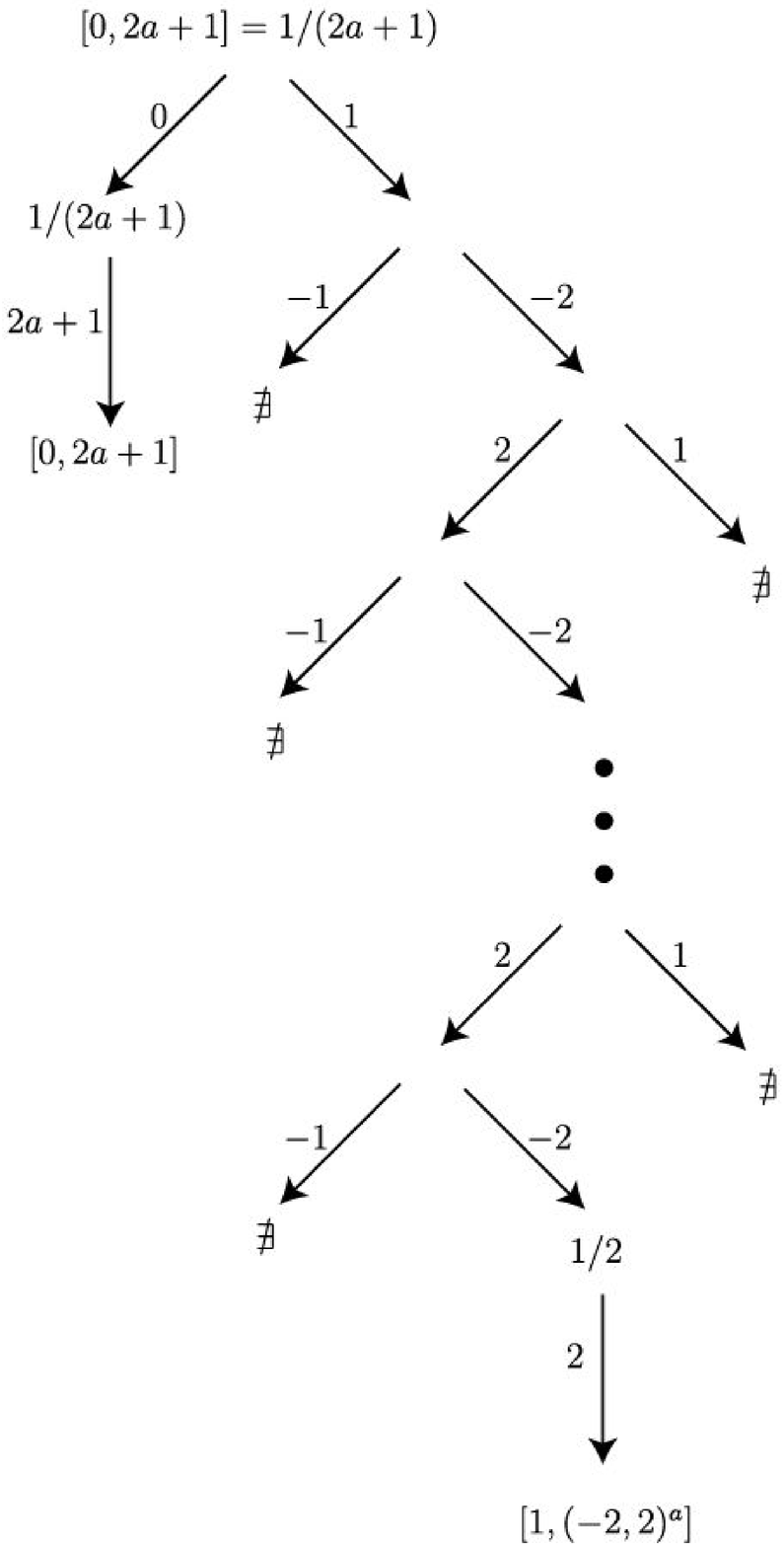}
\caption{\label{fig3}%
The binary tree for $[0,2a+1]$.}
\end{center}
\end{figure}

\underline{Subcase 3} ($a_0 = 2a+1$, $a \geq 1$): 
The binary tree is shown in Figure~\ref{fig3}.
 The two boundary slope continued fractions
$[0,a_0]$ and $[1,(-2,2)^a]$ are those given by substituting or not at position $0$.

\textbf{Case 2} ($n=1$):
Our goal is to show that the boundary slope continued fractions are among the fractions given
by substituting at position 0, at position 1, and by not substituting at all. The
result of substitution at position 0 will depend on whether $a_0$ is even or odd: 
\begin{eqnarray*}
{[}0,2a,a_1{]}    & \stackrel{\mbox{Sub.\ 1}}{\longrightarrow}  & {[}1, (-2,2)^{(a-1)},
-2, a_1 + 1 {]}
\\ {[}0,2a+1, a_1{]} & \stackrel{\mbox{Sub.\ 3}}{\longrightarrow} & {[}1, (-2,2)^a, -a_1-1
{]}
\end{eqnarray*}

Similarly, substitution at position 1 depends on the parity of $a_1$:
\begin{eqnarray*}
{[}0,a_0,2b{]}     & \stackrel{\mbox{Sub.\ 1}}{\longrightarrow} & {[}0, a_0+1,
(-2,2)^{(b-1)}, -2 {]} \\ 
{[}0,a_0, 2b+1{]} & \stackrel{\mbox{Sub.\ 3}}{\longrightarrow} &
{[}0, a_0+1, (-2,2)^b {]}
\end{eqnarray*}

\begin{figure}[h]
\begin{center}
\includegraphics[scale=0.43]{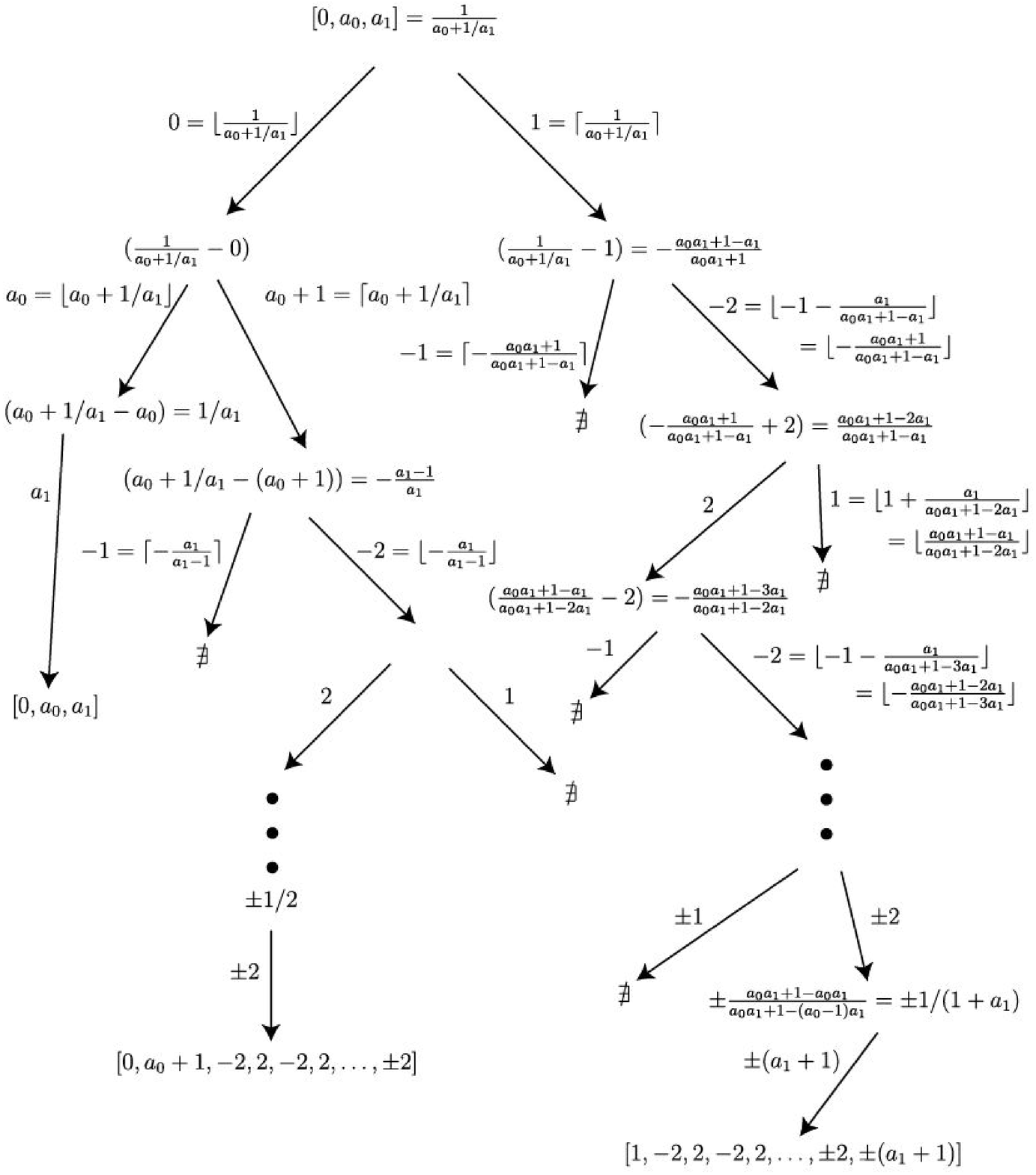}
\caption{\label{fig4}%
The binary tree for $[0,a_0,a_1]$.}
\end{center}
\end{figure}

As Figure~\ref{fig4} shows,
 these two boundary slopes, along with the original continued
fraction $[0, a_0, a_1]$ (no substitutions) are precisely those that arise in the
binary tree. Note that if, for example, $a_0$ or $a_1$ is $1$, then the $[0,a_0,a_1]$ leaf is
not in fact a boundary slope continued fraction. The point is that all leaves of the binary tree
are included in the set of continued fractions obtained by substitutions at non-adjacent 
positions. So, every boundary slope continued fraction appears in this set.

\textbf{Case 3} ($n=2$):
This case will illustrate how the induction works. There are five continued fractions given by 
substitutions at non-adjacent positions (compare with the example of Section~5.1): three
obtained by substitutions at positions 0, 1, and 2; one by substitutions at 0
and 2; and the original continued fraction itself (with no substitutions). Let us denote these
choices of substitutions by a sequence of three 0's and 1's where a 1 in the $i$th place 
denotes a substitution at that $i$th position. Thus, the five continued fractions will be
denoted $100$, $010$, $001$, $101$, and $000$. 

\begin{figure}[h]
\begin{center}
\includegraphics[scale=0.47]{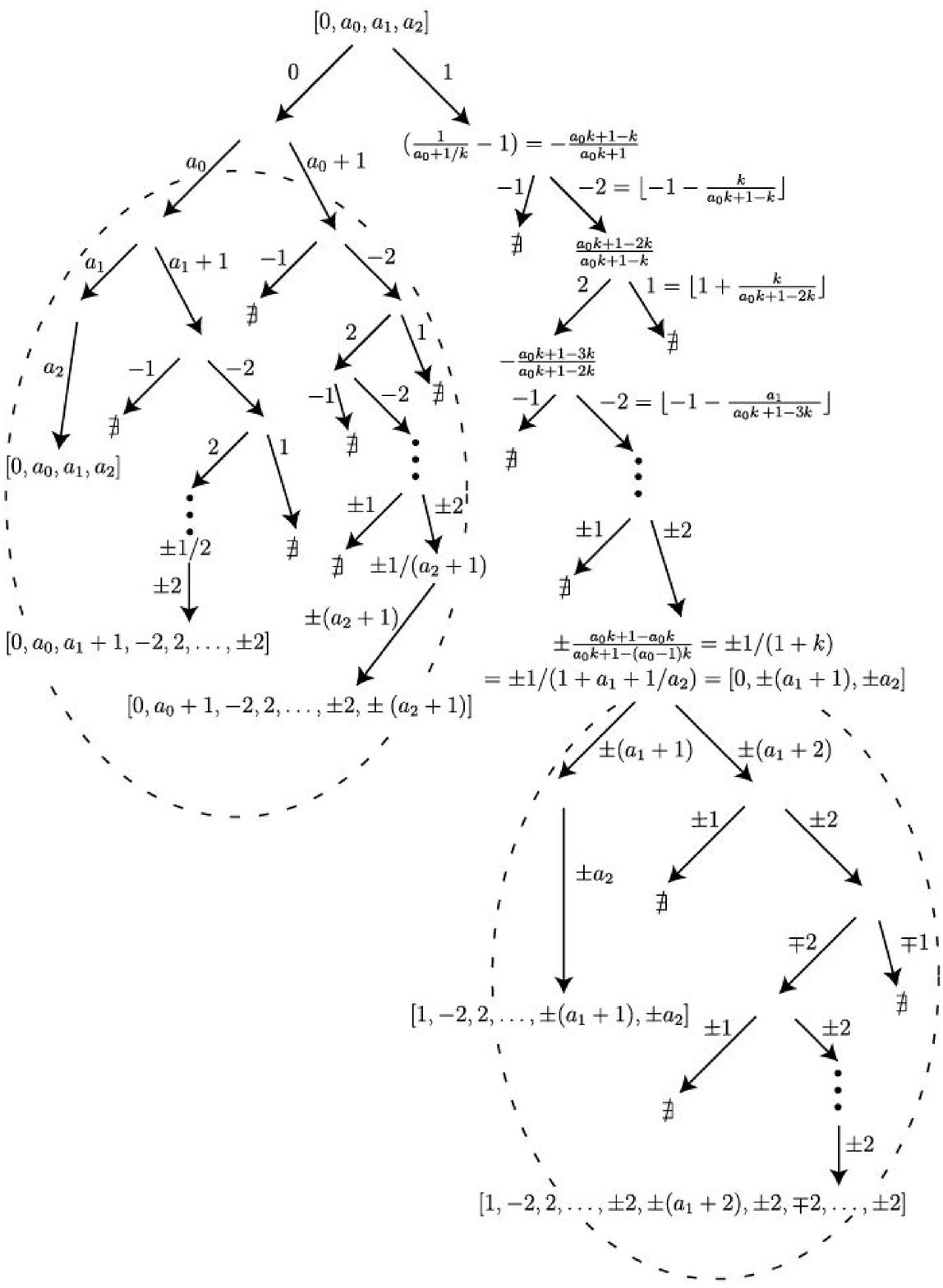}
\caption{\label{fig5}%
The $[0,a_0,a_1,a_2]$ tree is a union of two subtrees. 
}
\end{center}
\end{figure}

We can think of the 
binary tree (Figure~\ref{fig5}) 
as being a union of two subtrees. The one at left corresponds to making no
substitution at position 0. This subtree ends in the three boundary slopes which have:
no substitutions ($000$); substitution at position 1 ($010$); and substitution at position 2
($001$), i.e., the sequences that begin in $0$.
This subtree is essentially the same as that for the $[0,a_1,a_2]$ continued fraction
(compare Figure~\ref{fig4}) as we can obtain these three sequences by adding a $0$ at the
front of the three boundary slopes sequences $00$, $10$, and $01$ of that case. 
The other subtree corresponds to making a substitution at position 0 and no substitution at
position 1. This  subtree contains the remaining two boundary slopes: substitution at position
0 ($100$); and substitution at positions 0 and 2 ($101$), i.e.,
sequences that begin in $10$. This subtree is similar to that for $[0,a_2]$ (compare
Figure~\ref{fig2}) as it remains only to decide whether or not to substitute in the second 
position. Again, some of these five sequences may not result in a boundary slope continued
fraction, for example, if one of the $a_i$ is $1$. However, every leaf  of the tree will
be included in the set of continued fractions obtained by substituting at non-adjacent
positions.

\begin{figure}[h]
\begin{center}
\includegraphics[scale=0.26]{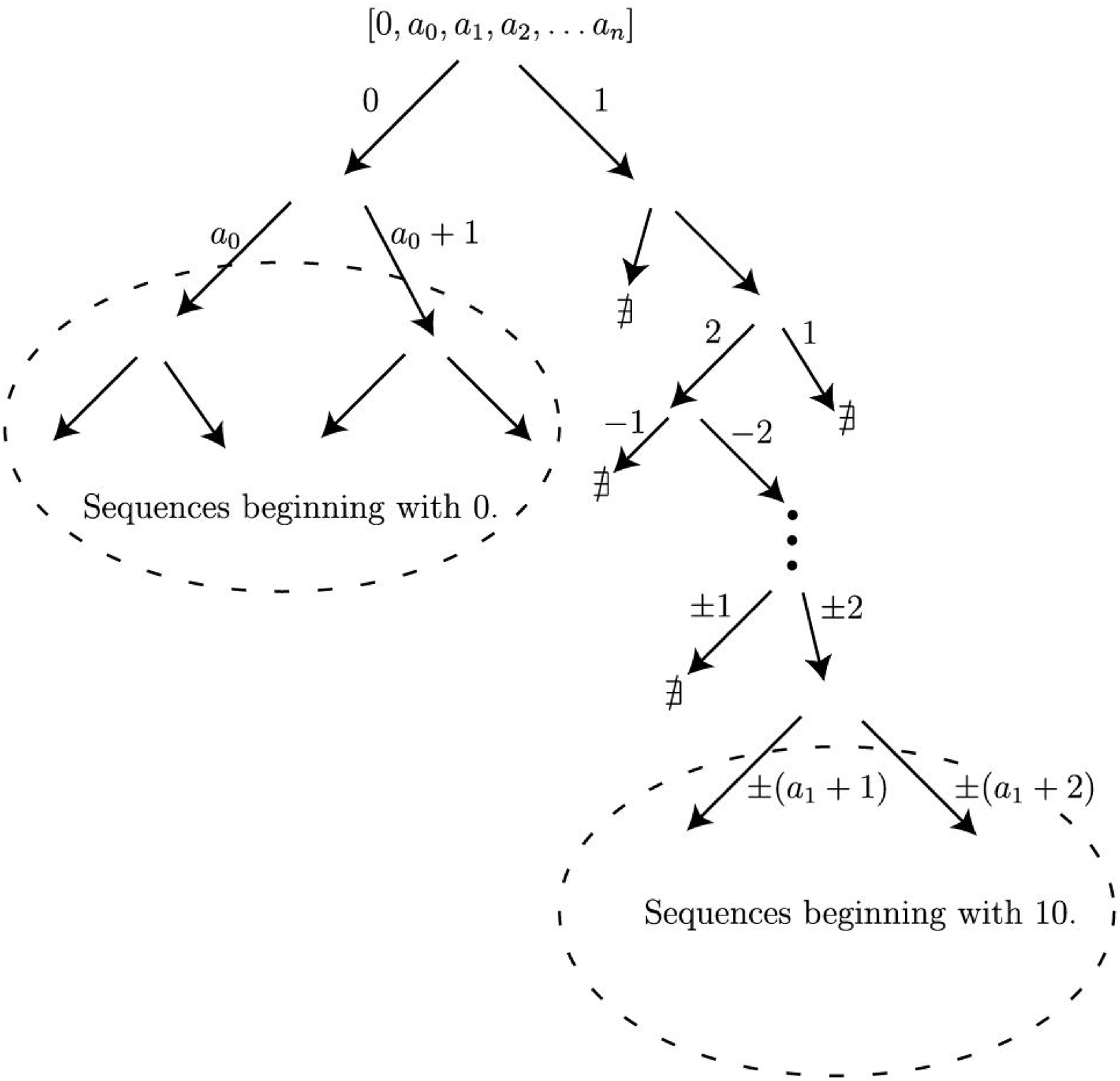}
\caption{\label{fig6}%
The general case also results in two subtrees.}
\end{center}
\end{figure}

\textbf{Case 4} ($n \geq 3$):
As in Case 3, we can decompose the binary tree (Figure~\ref{fig6})
 into two subtrees. One corresponds
to sequences that begin with $0$, the other to sequences beginning with $10$.  
The first will be, essentially, the tree that arises from the simple continued
fraction $[0,a_1,a_2, \ldots, a_n]$. By induction, the leaves of this subtree
correspond to non-adjacent substitutions in this simple continued fraction. 
By its placement in the $[0, a_0, a_1, \ldots, a_n]$ tree, this ensures that
the leaves of this part of the tree will correspond to continued fractions obtained
by substitution sequences into $[0,a_0, a_1, \ldots, a_n]$ that begin with $0$.

The other subtree is isomorphic to the tree that arises from the simple continued fraction
$[0,a_2, a_3, \ldots, a_n]$. By induction, the leaves of the subtree correspond to 
substitutions into this continued fraction. By its placement in the tree for
$[0,a_0,a_1,\ldots,a_n]$, the leaves here can be obtained by non-adjacent substitutions
into that continued fraction that begin with $10$.

Thus, every leaf of the binary tree and, therefore, every boundary slope continued fraction
can be obtained by non-adjacent substitutions into the simple continued fraction.
\QED

\setcounter{cor}{1}

\begin{cor} If $\frac p q =[0, a_0, a_1, \ldots, a_n]$ is a simple continued fraction, then $K(
p/q )$ has at most $F_{n+2}$ boundary slopes where $F_n$ is the $n$th Fibonacci number.
\end{cor}

\proof
We have shown that the boundary slope continued fractions lie among those given by substitution
at non-adjacent positions which in turn are in bijection with sequences of $n+1$ 0's or 1's
containing  no pair of consecutive 1's. Thus the number of boundary slopes is at most $P_n$,
where $P_n$ is the number of $0$, $1$ sequences of length $n+1$ with no consecutive 1's. We will
show that
$P_n = F_{n+2}$ by induction.

There are two base cases. If $n=0$,
there are two sequences: $0$ and $1$. So, $P_0 = 2 = F_2$.
For $n=1$, there are three sequences: $00$, $10$, and $01$. So, $P_1 = 3 = F_3$.

For the inductive step, sequences of length $n+1$ are obtained by either adding a $0$ to
the beginning of a $n$ sequence or $10$ to the beginning of a $n-1$ sequence. 
Thus $P_n = P_{n-1} + P_{n-2} = F_{n+1} + F_{n} = F_{n+2}$. \QED

In general, $F_{n+2}$ is an overestimate since the continued fractions obtained by substitutions
will not necessarily have terms at least two in absolute value. In particular, if the
simple continued fraction includes any 1's, then the continued fraction obtained by making
no substitutions ($000 \dots 0$) will not be a boundary slope continued fraction. Moreover,
different boundary slope continued fractions could result in the same boundary slope.
For example, this will occur when, in the simple continued fraction, we have two equal
terms separated by an even distance:
$a_i = a_{i+2k}$.

\section{Maximum and Minimum Boundary Slopes}
	
In this section, we will show how one can calculate the maximum and minimum boundary slopes.
We will refer to the minimum value of $b^+-b^-$ as $b_1$ and the maximum as $b_2$. 
Further, we will refer to the components of $b_1$ as $b_1^+$ and $b_1^-$.  Similarly for $b_2$.
	
The key observation is that, if we begin with a simple continued fraction, applying a
substitution at an even position will decrease $b^+ - b^-$ and, hence, the boundary slope, while
applying a substitution at an odd position will increase the boundary slope. Thus, we can
minimize the boundary slope by substituting at each even position (and no odd positions). Note
that this will result in a  continued fraction where each term is at least two in absolute value.
Indeed, the even position terms of the original simple continued fraction will be replaced by a
sequence of
$\pm 2$'s while the terms in the odd positions will be augmented in absolute value by at least
one. 
	
However, we need a way to count the resulting $b_1^+$ and $b_1^-$ when we make substitutions at
each even position. We will do this by focusing on the clusters of $(\pm 2, \mp 2)$,
examining what occurs near them individually, and then summing up the results.  
We replace an even number, $2k$, with $(\pm 2, \mp
2)^{|k|-1}, \pm 2$, or $2|k|-1$ terms.  
For a positive odd number
$2k+1$, we replace it with $(\pm 2, \mp 2)^k$
or $2k$ terms and a negative odd $-(2k+1)$ is also replaced with  $2k$ terms.
We can combine these cases by observing that a term $a_i$ is replaced with $|a_i|-1$ terms.
If we think of making the substitutions at position $0$, $2$, $4, \ldots$ in turn, then,
each substitution will not affect the magnitude of later $a_{2i}$ but may change their signs.
Thus, if we begin with the simple continued fraction $[0,a_0, \ldots, a_n]$, each $a_{2i}$
will be replaced by $a_{2i}-1$ terms.
	
All we need now, then, is to count the number of terms that appear between the strings of
$\pm2$'s.  There will be one such term between every string of $\pm2$'s, and possibly one at the
tail end of the continued fraction.  Specifically, for the simple continued fraction $[0, a_0,
\dots, a_n]$ where $n$ is even, there are 
$\frac n 2$ odd position terms, $a_1, a_3, \ldots, a_{n-1}$.  Similarly, if $n$ is odd, then we
must have
$\frac{n+1}{2}$ odd terms.  More concisely, there are $\ceil{\frac{n}{2}}$ odd terms. 
Now, we have all the ingredients necessary to calculate $b_1^-$.
	
	\begin{eqnarray*}
		b_1^-	&=&	\ceil{\frac{n}{2}} + \sum_{i=0}^{\floor{\frac{n}{2}}}  \left(a_{2i}-1\right) \\
				&=&	\ceil{\frac{n}{2}} - \left(\floor{\frac{n}{2}}+1\right) + \sum_{i=0}^{\floor{\frac{n}{2}}} a_{2i}  \\
				&=&	\sum_{i=0}^{\floor{\frac{n}{2}}} a_{2i} \mbox{\qquad if $n$ is odd} \\
				&=&	-1 + \sum_{i=0}^{\floor{\frac{n}{2}}} a_{2i} \mbox{\qquad if $n$ is even}
	\end{eqnarray*}
	
So, surprisingly, $b_1^-$ is equal to the sum of the even terms, less one if $n$ is even. 
Recall that, the simple continued fraction is $[0, a_0, \dots,
a_n]$ so that $n$ is one more than the number of partial quotients. Note that $b_1^+ = 0$
as all partial quotients are now matched to the $[-+-+ \cdots]$ pattern.

Similarly, the maximum boundary slope is given by substituting at each odd position,
and we can use the exact same logic to find $b_2^+$. In fact, the formula itself is nearly
identical.  The sole differences are the terms we sum over, and the fact that we don't add
$\floor{\frac{n}{2}}$, but instead $\ceil{\frac{n+1}{2}}$.  We then get:
	
	\begin{eqnarray*}
		b_2^+	&=&	\ceil{\frac{n+1}{2}} + \sum_{i=0}^{\floor{\frac{n-1}{2}}}  \left(a_{2i+1}-1\right) \\
				&=&	\ceil{\frac{n+1}{2}} - \left(\floor{\frac{n-1}{2}}+1\right) + \sum_{i=0}^{\floor{\frac{n-1}{2}}} a_{2i+1}\\
				&=&	\sum_{i=0}^{\floor{\frac{n-1}{2}}} a_{2i+1} \mbox{\qquad if $n$ is odd} \\
				&=&	1 + \sum_{i=0}^{\floor{\frac{n-1}{2}}} a_{2i+1} \mbox{\qquad if $n$ is even}
	\end{eqnarray*}
	
Now that we have $b_1$ (which is $-b_1^-$, since $b_1^+=0$) and $b_2$ (which is simply $b_2^+$),
we can calculate the maximum and minimum boundary slopes in terms of $b_0^+$ and $b_0^-$.  The
minimum boundary slope is $2\big((0 - b_1^-) - (b_0^+ - b_0^-)\big) = -2b_1^- - 2(b_0^+ -
b_0^-)$.  Similarly, the maximum is $2\big((b_2^+ - 0) - (b_0^+ - b_0^-)\big) = 2b_2^+ - 2(b_0^+
- b_0^-)$.
		
	\section{Proof of Theorem~\ref{thmain}}
	
	In this section we prove our main theorem, that twice the crossing number
of a $2$--bridge knot $K$ is equal to the diameter of the boundary slopes.
	
\setcounter{theorem}{0}

	\theorem{For $K$ a $2$--bridge knot, $D(K) = 2c(K)$.}
	
\proof Firstly, calculating the crossing number of a rational knot is simple.  If $[0, a_0,
\dots, a_n] = \frac p q$ is the simple continued fraction for $K\left(\frac p q\right)$, then
$c\left(K\left(\frac p q\right)\right) = \sum_{i=0}^n a_i$.
	
The diameter of $B(K)$ is also easy to calculate.  If we use the $b_1$ and $b_2$ from the
previous section, we get $D(K) = 2b_2^+ - 2(b_0^+ - b_0^-) - \big( -2b_1^- - 2(b_0^+ - b_0^-)
\big) = 2b_2^+ + 2b_1^-$.  At this point, $b_1^-$ and $b_2^+$ may vary depending on whether $n$
is even or odd.  However, the differences cancel each other out in either instance, leaving us
with
	\begin{eqnarray*}
		D(K)	&=&	2\sum_{i=0}^{\floor{\frac{n-1}{2}}} a_{2i+1} + 2\sum_{i=0}^{\floor{\frac{n}{2}}}
a_{2i} \\
				&=&	2\sum_{i=0}^n a_i
	\end{eqnarray*}
	
	This concludes the proof that $2c(K) = D(K)$.  \QED
	
	


\begin{thebibliography}{99}
	
	\bibitem{A} C.C.\ Adams. \textit{The Knot Book: An Elementary
Introduction to the Mathematical Theory of Knots}. American Mathematical
Society (2004).
	
	\bibitem{C} P.\ Cromwell. \textit{Knots and Links}. Cambridge University
Press (2004).
	
	\bibitem{ES} C.\ Ernst and D.W.\ Sumners. The growth of the number of
prime knots. \textit{Math.\ Proc.\ Cambridge Phil.\ Soc.} \textbf{102}
(1987) 303--315.
	
	\bibitem{HT} A.\ Hatcher and W.\ Thurston. Incompressible surfaces in
2-bridge knot complements. \textit{Inventiones Math.} \textbf{79} (1985)
225--246.
	
	\bibitem{I} K.\ Ichihara, Private Communication.

 \bibitem{IM} K.\ Ichihara and S.\ Mizushima. 
	Crossing number and diameter of boundary slope set of Montesinos knot. Preprint (2005). 
    math.GT/0510370
	
\bibitem{IS} M.\ Ishikawa and K.\ Shimokawa. Boundary slopes and crossing
numbers of knots. Preprint (2002).
	
	\bibitem{L} R.\ Langford. M.S. Thesis, CSU, Chico (2005).
	
\bibitem{Mo}  L.\ Moser. Elementary surgery along a torus knot. 
\textit{Pacific J.\ Math.}  \textbf{38} (1971) 737--745.

\bibitem{Mu} K.\ Murasugi. 
 On the braid index of alternating links.  
\textit{Trans.\ Amer.\ Math.\ Soc.}  \textbf{326}  (1991) 237--260.
	
\end{thebibliography}
\end{document}